\documentclass[11pt]{article}
\usepackage{amsthm,amssymb,amscd,amsmath,latexsym,indentfirst,color,amsfonts}
\usepackage[mathcal]{eucal}

\theoremstyle{plain}
\newtheorem{theorem}{Theorem}[section]
\newtheorem{prop}[theorem]{Proposition}
\newtheorem{coro}[theorem]{Corollary}
\newtheorem{lem}[theorem]{Lemma}
\newtheorem{defi}[theorem]{Definition}

\newcommand{\ble}{\begin {lem}}
\newcommand{\ele}{\end {lem}}
\newcommand{\bpr}{\begin {prop}}
\newcommand{\epr}{\end {prop}}
\newcommand{\bthm}{\begin {theorem}}
\newcommand{\ethm}{\end {theorem}}
\newcommand{\bco}{\begin {coro}}
\newcommand{\eco}{\end {coro}}
\newtheorem{rem}[theorem]{Remark}

 \DeclareMathOperator\Aut{Aut}
\DeclareMathOperator\Hom{Hom}

\DeclareMathOperator\an{and}
\DeclareMathOperator\Stab{Stab}\DeclareMathOperator\Rad{Rad}
\DeclareMathOperator\Supp{Supp}

\def\sl{\frak s\frak l}   

\def\bt{\beta}
\def\si{\sigma}
\def\al{\alpha}\def\Ga{\Gamma}
\def\w{\wedge}\def\rt{\rightarrow}\def\st{\subset}\def\ot{\otimes}\def\op{\oplus}
\def\wh{\widehat}

\newcommand{\be}{\begin {equation}}
\newcommand{\ee}{\end {equation}}
\newcommand{\bp}{\begin {proof}}
\newcommand{\ep}{\end {proof}}
\newcommand{\bee}{\begin {equation*}}
\newcommand{\eee}{\end {equation*}}
\newcommand{\lb}{\label}

\def\sl{\frak s\frak l}

\begin{document}

\title{Multiplicity free gradings on semisimple Lie and Jordan algebras and skew root systems}

\author{Gang Han \\ Department of Mathematics \\ Zhejiang University\\China \thanks{Corresponding author. Email address:mathhgg@zju.edu.cn(G.Han). His work is supported by Zhejiang Province Science Foundation, grant No. LY14A010018.}
\and Yucheng Liu \\ Department of Mathematics \\ Zhejiang University\\China
\and Kang Lu \\ Department of Mathematics \\ Zhejiang University\\China
 }

\date{Nov. 11, 2016}

%\dedicatory{This paper is dedicated to our authors.}
%\thanks{{\em 2000 Mathematics Subject Classification.} Primary 17B20, 17B10.}

\maketitle

{\small \noindent \textbf{Abstract}.

 A $G$-grading on an algebra is called multiplicity free if each homogeneous component of the grading is 1-dimensional, where $G$ is an abelian group. We introduce skew root systems of Lie type and skew root systems of Jordan type respectively, and use them to construct multiplicity free gradings on semisimple Lie algebras and on semisimple Jordan algebras respectively. Under certain conditions the corresponding Lie (resp. Jordan) algebras are simple. Three families of skew root systems of Lie type (resp. of Jordan type) are constructed and the corresponding Lie (resp. Jordan) algebras are identified. This is a new approach to study abelian group gradings on Lie and Jordan algebras.

}


\tableofcontents
%\noindent \textbf{2000 Mathematics Subject Classification}:17B10\bigskip

%\noindent \textbf{Key words}:  }
\section{Introduction}
 \setcounter{equation}{0}\setcounter{theorem}{0}

People are interested in  gradings on nonassociative algebras in the last years. They have done a lot of work in this area and find some important applications. An algebra considered in this paper will be either an associative algebra,  a Jordan algebra or a Lie algebra over some algebraically closed field $F$ of characteristic 0, and will always be assumed to be finite dimensional over $F$ unless otherwise specified.  An associative or Jordan algebra will always be assumed to be unital, and a grading on some algebra $A$ will always refer to an abelian group grading, i.e., a grading on $A$ by some abelian group. We will focus on gradings on semisimple Lie and Jordan algebras. Good surveys for such gradings can be found in Chapter 3.3 of \cite{ov}, \cite{k} and \cite{ek2}.

The classical Cartan decomposition of a semisimple Lie algebra $L$ is in fact a grading on $L$ by a free abelian group of the same rank as that of $L$. Alekseevskii found several classes of interesting gradings on simple Lie algebras by elementary abelian groups and computed the corresponding Weyl groups in \cite{a}. Since late 1980's, Patera and his collaborators began to study abelian group gradings on Lie algebras systematically in \cite{pz}, \cite{hpp} and etc, and they classified most of the fine gradings on classical simple Lie algebras. Later in \cite{bsz} Bahturin and his collaborators described all the abelian group gradings on classical simple Lie algebras except the simple Lie algebra of type $\textbf{D}_4$ and on several classes of simple Jordan algebras as well. In \cite{bk} Bahturin and Kochetov classified the isomorphism classes of abelian group gradings on all the simple Lie algebras except the simple Lie algebra of type $\textbf{D}_4$ in terms of numerical and group-theoretical invariants. Draper and his collaborators classified all the fine gradings on the simple Lie algebra of type $\textbf{D}_4$ in \cite{dv}, classified all the gradings on the simple Lie algebra of type $\textbf{G}_2$ in \cite{dm1} and classified all the fine gradings on the Albert algebra and on the simple Lie algebra of type $\textbf{F}_4$ in \cite{dm2}. Notice that many of the above mentioned works make use of knowledge on gradings of associative algebras.\bigskip

There is a duality between gradings on an algebra $A$ and quasitori in $\Aut(A)$, which will be reviewed in next section. Isomorphic gradings on $A$ correspond to conjugate quasitori in $\Aut(A)$, and fine gradings on $A$ correspond to maximal quasitori in $\Aut(A)$. The Weyl group of a grading is defined to describe its symmetry, and in \cite{ek} Elduque and Kochetov   computed the Weyl groups for all the fine gradings on all the classical simple Lie algebras except $\textbf{D}_4$.  Also see Section 2 of \cite{dm2} for some details.

Let $A$ be an algebra. A grading on it is called \textit{multiplicity free} if each component of the grading is 1-dimensional. If  $A$ is an associative algebra or a Jordan algebra, then a grading on it is called \textit{division} if each nonzero homogeneous element of $A$ is invertible. An associative or Jordan algebra with a multiplicity free division $G$-grading is called an associative or Jordan $G$-torus respectively in \cite{ny}. A twisted group algebra has a canonical multiplicity free division grading. Conversely, it can be shown easily that an associative algebra with a multiplicity free division grading must be a twisted group algebra. See Proposition \ref{3a}. A twisted group algebra over an abelian group $G$ corresponds to  an alternating bicharacter $\bt$ on $G$. Such $(G,\bt)$ is called a symplectic abelian group. There is a one-to-one correspondence between twisted group algebras over the abelian group $G$ and alternating bicharacters on $G$. The Weyl group of a twisted group algebra is studied in \cite{op}. See Section 3 for details.

If  $A$ is a Lie algebra, we call a grading on it \textit{semisimple} if each  homogeneous element of $A$ is semisimple. The classical Cartan decomposition of a semisimple Lie algebra is not a semisimple grading as each root space consists of nilpotent elements, but the grading on a semisimple Lie algebra induced by the action of a finite maximal quasitorus in $\Aut(A)$ is a  semisimple grading.

 A root system $R$ is a finite set of elements in a Euclidean space $V$ satisfying certain axioms. There is a one-to-one correspondence between reduced and irreducible root systems and simple Lie algebras over algebraic fields of characteristic zero (up to isomorphism). In particular, given a root system one will construct a semisimple Lie algebra. Analogously, we define skew root systems of Lie type and skew root systems of Jordan type respectively, and use them to construct multiplicity free semisimple gradings on semisimple Lie algebras and multiplicity free division gradings on semisimple Jordan algebras respectively. These are the main results of this paper, which are contained in Section 4 and 5.

  Given a symplectic abelian group $(G,\bt)$, one can construct a twisted group algebra $A$ over $G$. One knows that $A$ has a canonical multiplicity free division grading. The associative algebra $A$ can be regarded as a Lie algebra, denoted by $A^-$, with the Lie bracket $[a,b]=ab-ba$; and similarly $A$ can be regarded as a Jordan algebra, denoted by $A^+$, with the Jordan product $a\circ b=\frac{1}{2}(ab+ba)$. Since $G$ is abelian, the grading on $A$ is compatible with the multiplication in $A^-$ or $A^+$. Thus $A^-$ and $A^+$ also have multiplicity free gradings induced from that on $A$. A set $R$ of elements in the symplectic abelian group $(G,\bt)$ satisfying certain axioms is called a skew root system of Lie type. Given such a $(G,\bt,R)$, one always obtains a Lie algebra $L(R)$ with a multiplicity free semisimple grading such that $L(R)$ is a graded Lie subalgebra of $A^-$. Similarly,  a set $R$ of elements in the symplectic abelian group $(G,\bt)$ satisfying another set of axioms is called a skew root system of Jordan type. Given such a $(G,\bt,R)$, one always obtains a  Jordan algebra $J(R)$ with a multiplicity free division grading such that $J(R)$ is a graded Jordan subalgebra of $A^+$. If $G$ is a finite abelian group, then it is shown in Proposition \ref{413} and \ref{414} that the algebras $L(R)$ and $J(R)$ are semisimple. Moreover, it is shown in Theorem \ref{411} that if $(G,\bt,R)$ is a reduced and indecomposable skew root system of Lie type, then $L(R)$ is always simple; and if $(G,\bt,R)$ is a skew root system of Jordan type, then it is proved in Proposition \ref{412} that $J(R)$ is always graded-simple. So one can construct (maybe infinite dimensional) simple Lie algebras (resp. graded-simple Jordan algebras) by means of skew root systems of Lie type (resp. of Jordan type). See Section 4 for details.

  Finally in Section 5 we constructed 3 families of skew root systems of Lie type (resp. of Jordan type) in finite symplectic abelian groups, and constructed the corresponding graded Lie algebras (resp. graded Jordan algebras). Most of these algebras are simple. Although such gradings are already known before, they have never been constructed by means of skew root systems.  As far as we know, all the multiplicity free semisimple gradings on simple Lie algebras and multiplicity free division gradings on special simple Jordan algebras are included in them. In particular, exceptional simple Lie algebras do not have multiplicity free semisimple gradings. So we believe that all the multiplicity free semisimple gradings on simple Lie algebras and all the multiplicity free division gradings on special simple Jordan algebras can be constructed by means of skew root system of Lie type and of Jordan type respectively.

\section{Some preliminaries on gradings}
 \setcounter{equation}{0}\setcounter{theorem}{0}

 \noindent \textbf{2.1.} We will first briefly review gradings on a nonassociative algebra by abelian groups. Let $F$ be an algebraically closed field of characteristic zero. All the algebras in the paper are assumed to be  finite dimensional over $F$ unless otherwise specified.

Let $A$ be a algebra and $G$ be a finitely generated additive abelian group. A $G$-\textit{grading} $\Gamma$ on $A$ is the decomposition of
$A$ into a direct sum of subspaces $$\Gamma: A=\oplus_{g\in G}
\ A_g$$ such that
$$A_g\cdot A_h\st A_{g+h},\ \forall\ g,h\in G.$$ We also say that $A$ is $G$-\textit{graded}. The subset
$R=\{g\in G|A_g\neq 0\}$ of $G$ is called the \textit{support} of this grading and is denoted $\Supp\ \Ga$. For any $g\in R$, $A_g$ is called the \textit{homogeneous component} of degree $g$ and each element in $A_g$ is a \textit{homogeneous} element. We will always
assume that $G$ is generated by $R$, otherwise it could be
replaced by its subgroup generated by $R$. So $G$ is always a finitely generated abelian group. Denote this grading by $(\Ga,G)$, or simply by $\Ga$. All the gradings mentioned in the paper will be some $G$-grading with $G$ a (finitely generated) additive abelian group.

Let $\Ga$ be a $G$-grading on the algebra $A$. Let $\Aut(A)$  be the algebraic group of automorphisms of $A$. The group $\Aut(\Ga)$ consists of all the automorphisms in $\Aut(A)$ that permute the homogeneous components of $\Ga$. The group $\Stab(\Ga)$ consists of all the automorphisms in $\Aut(A)$ that stabilize each homogeneous component of $\Ga$.  Clearly $\Stab(\Ga)$ is a normal subgroup of $\Aut(\Ga)$, and let $$W(\Ga)=_{def}\Aut(\Ga)/\Stab(\Ga)$$ be the \textit{Weyl} group of the grading $\Ga$, which describes the symmetry of the grading. See also \cite{ek}.

\noindent \textbf{2.2.} Now we recall the duality between abelian group gradings and abelian group actions.

An abelian linear algebraic group consisting of semisimple elements is called a quasitorus. Recall that $F$ is an algebraically closed field of characteristic zero. Assume $T$ is a quasitorus over $F$. Let $\wh{T}$ be the set of algebraic group homomorphisms from $T$ to $F^\times$. Then $\wh{T}$ is an abelian group, called the \textit{character group} or \textit{dual group} of $T$ with addition defined by
\[ (\al+\bt)(a)=\al(a)\cdot\bt(a),\ for\ any\ \al,\bt\in \wh{T}; \ for\ any\ a\in T. \]
One knows that $T\mapsto \wh{T}$ is a contravariant functor from the category of quasitori over $F$ to the category of finitely generated abelian groups and is a duality  between these two categories. The inverse functor maps a finitely generated abelian group $G$ to the set $\wh{G}=Hom(G, F^\times)$ of group homomorphisms, which has a canonical structure of quasitorus over $F$. This is a version of Pontryajin duality in the setting of algebraic groups.

Recall that $A$ is a finite dimensional non-associative algebra. The group $\Aut(A)$ is a linear algebraic group over $F$. Assume $T\st\Aut(A)$ is a quasitorus. Then the action of $T$ on $A$ induces a $\wh{T}$-grading on $A$:
\be \Ga:A=\op_{g\in \wh{T}}A_g. \ee  If $T^{'}$ is another quasitorus in $\Aut(A)$, then $T$ and $T^{'}$ induces isomorphic group gradings on $A$ if and only if $T$ and $T^{'}$ are conjugate in $\Aut(A)$.

Conversely given a grading $A=\oplus_{g\in G}\ A_g$ with $G$ a finitely generated abelian group, let $\wh{G}$ be the dual group of $G$, which is a  (multiplicative) abelian group consisting of homomorphisms from $G$ to $F^\times$. Then the grading induces a $\wh{G}$-action on $A$:
\be \si\cdot X=\si(g)X,\ for\ all\ \si\in G, X \in A_g.\ee This action gives a homomorphism $\wh{G}\rt \Aut(A)$, which is injective as the support of the grading generates $G$. So the homomorphism embeds $\wh{G}$ as a quasitorus in $\Aut(A)$ and we will identify $\wh{G}$ with its image in $\Aut(A)$.

So one has the one-to-one correspondence between abelian group gradings on $A$ and quasitori in $\Aut(A)$. Two gradings on $A$ are isomorphic if and only if the respective quasitori are conjugate in $\Aut(A)$. A $G$-grading on $A$ is called fine if $\widehat{G}$ is a maximal quasitorus in $\Aut(A)$. See Section 2 of \cite{dm2}.

For any quasitorus $T\st \Aut(A)$, let $N(T)$ and $Z(T)$ be respectively the normalizer and centralizer of $T$ in $\Aut(A)$. Define \[W(T)=N(T)/Z(T)\] and call it the Weyl group of $T$ with respect to $ \Aut(A)$, which is always a finite group.

Assume that $A$ is a semisimple (associative, Lie or Jordan) algebra. If  $T\st \Aut(A)$ is a maximal quasitorus, then $Z(T)=T$ and the action of $T$  induces a fine $G$-grading $\Ga$ on $A$ with $G=\wh{T}$. Then it is not hard to verify that $N(T)=\Aut(\Ga)$, $Z(T)=\Stab(\Ga)$ and thus \be\lb{l1}W(T)=W(\Ga).\ee A proof in the case that $A$ is a simple Lie algebra can be found in Proposition 2.4 of \cite{h1}, which applies to the general case as well.\bigskip

\section{An isomorphism between $H^2(G,F^\times)$ and $\w^2(G,F^\times)$}
 \setcounter{equation}{0}\setcounter{theorem}{0}
Recall that $F$ is an algebraically closed field with characteristic 0, and $G$ is a finitely generated abelian group.

It is well known that the second integral homology group of $G$ is isomorphic to the second exterior power $ \w^2 G$ of
$G$ with $G$ viewed as a $\mathbb{Z}$-module, specifically, \be\lb{a} H_2(G,\mathbb{Z})\rt \w^2 G, [g,h]\mapsto g\w h\ee is the isomorphism. By Universal Coefficient Theorem (cf. Exercise 6.1.5 of \cite{w}) one has a split exact sequence  \be\lb{b} \CD 1\rt \textrm{Ext}_\mathbb{Z}^1(H_1(G,\mathbb{Z}),F^\times)\rt H^2(G,F^\times)@>\psi>> Hom_\mathbb{Z}(H_2(G,\mathbb{Z}),F^\times)\rt 1.\endCD\ee

Let $\w^2 (G,F^\times)$ denote the abelian group of alternating $\mathbb{Z}$-bilinear maps $\bt:G\times G\rt F^\times$, i.e.,  $\bt$ is $\mathbb{Z}$-bilinear and $\bt(g,g)=1$ for any $g\in G$. An element in $\w^2 (G,F^\times)$ will be called an \textit{alternating bicharacter} of $G$. By (\ref{a}) one has \[Hom_\mathbb{Z}(H_2(G,\mathbb{Z}),F^\times)\cong \w^2 (G,F^\times).\] As $G$ is abelian, $H_1(G,\mathbb{Z})\cong G$. Thus (\ref{b}) becomes  \be \CD 1\rt \textrm{Ext}_\mathbb{Z}^1(G,F^\times)\rt H^2(G,F^\times)@>\psi>> \w^2 (G,F^\times)\rt 1.\endCD\ee

As $F$ is algebraically closed,  $F^\times$ is divisible and $\textrm{Ext}_\mathbb{Z}^1(G,F^\times)=0$. So one has
 \ble
The map $\psi$ is a $\mathbb{Z}$-linear isomorphism.
 \ele

 As $H^2(G,F^\times)\cong Z^2(G,F^\times)/B^2(G,F^\times)$, the isomorphism $\psi$ can be written as a short exact sequence of abelian groups

 \be\lb{c}\CD  1\rt B^2(G,F^\times)\rt Z^2(G,F^\times)@>\Psi>>  \w^2(G,F^\times)\rt 1,\endCD\ee where the map $\Psi$ is \be \lb{}\Psi: Z^2(G,F^\times)\rt \wedge^2(G,F^\times), \Psi(\xi)(a, b)= \xi(a,b) \xi(b,a)^{-1}.\ee

If $\bt\in  \wedge^2(G,F^\times)$, $\xi\in Z^2(G,F^\times)$ and $\Psi(\xi)=\bt$, then we refer to $\xi$ as a 2-cocycle corresponding to $\bt$.

 As $G$ is abelian, elements in $L^2(G,F^\times)=\Hom_\mathbb{Z}(G\ot G,F^\times)$ are $\mathbb{Z}$-bilinear maps $G\times G\rt F^\times$ (called bicharacters of $G$), which are clearly 2-cocycles in $Z^2(G, F^\times)$. Then $\Psi$ in (\ref{c}) restricts to a  homomorphism
\be \Psi:L^2(G,F^\times)\rt  \w^2(G,F^\times),\ee which is surjective as $G$ is finitely generated. Denote the kernel of $\Psi$  by $ \textrm{Sym}^2(G,F^\times)$.
\ble\lb{pd}
One has a short exact sequence of abelian groups
\be  1\rt \textrm{Sym}^2(G,F^\times)\rt L^2(G,F^\times)\rt \wedge^2(G,F^\times)\rt 1.\ee

\ele

Now we recall the definition of twisted group algebras. For any 2-cocycle $\xi\in Z^2(G,F^\times)$, let $F^\xi G$ be the $F$-vector space with basis $\{u_g|g\in G\}$ and the multiplication on the basis elements is defined by \be u_a \cdot u_b=\xi(a,b) u_{a+b}, for\ all\ a,b\in G.\ee Extending the multiplication to $F^\xi G$ bilinearly. Then $F^\xi G$ is called a twisted group algebra of $G$, which is an associative algebra with a canonical multiplicity free division $G$-grading: \[F^\xi G=\bigoplus_{g\in G} F\cdot  u_g. \]

One knows that cohomologous 2-cocycles in $ Z^2(G,F^\times)$ define isomorphic $G$-graded  associative algebras, and the isomorphic classes of twisted group algebras of $G$ over $F$ are in one-to-one  correspondence with $H^2(G,F^\times)$. The above results can be found in Section 1, Chapter 2 of \cite{ka2}.

We have showed that $H^2(G,F^\times)\cong \wedge^2(G,F^\times)$ thus isomorphic classes of twisted group algebras of an abelian group $G$ are in one-to-one correspondence with $\wedge^2(G,F^\times)$.

 Let $G$ be a finitely generated (additive) abelian group with an alternating bicharacter $$\bt:G\times G\rt F^\times.$$ Then $(G,\bt)$ is called a \textit{symplectic abelian group}. Let \[\Rad(\bt)=\{g\in G|\bt(g,h)=1, \ for\ all\ h\in G \}\] be the \textit{radical} of $(G,\bt)$, which is a subgroup of $G$. We say that $(G,\bt)$ is \textit{nonsingular} if $\Rad(\bt)=0$. For any two subsets $A,B\subseteq G$, we say that $A$ and $B$ are \textit{orthogonal} and denoted $A\perp B$, if $\bt(a,b)=1$ for any $a\in A, b\in B$.

  Given $(G,\bt)$, there corresponds to it a twisted group algebra $F^\xi G$, where $\bt=\Psi(\xi)$. Note that if $\xi^{'}\in Z^2(G,F^\times)$ is another 2-cocycle with $\bt=\Psi(\xi^{'})$ then $F^{\xi^{'}} G$ and $F^\xi G$ are $G$-graded isomorphic. If $G$ is finite then $F^\xi G$ is semisimple, and $F^\xi G$ is simple if and only if $(G,\bt)$ is nonsingular. See Lemma 2.7 and Corollary 2.10 in Chapter 8 of \cite{ka3} for the proof.

 Let $(G_i,\bt_i),\ i=1,2$ be two symplectic abelian groups.  Let $G=G_1\op G_2$ and define the bicharacter $\bt$ on $G$ by \be\lb{}\bt((a_1,a_2),(b_1,b_2))=\bt_1(a_1,b_1) \bt_2(a_2,b_2)\ee with $a_i,b_i\in G_i$. Then $(G,\bt)$ is called the\textit{ orthogonal direct sum} of $(G_i,\bt_i),\ i=1,2$. It is clear that $\Rad(\bt)=\Rad(\bt_1)\op \Rad(\bt_2)$ thus $\Rad(\bt)=0$ if and only if $\Rad(\bt_i)=0$ for $i=1,2$. So $(G,\bt)$ is nonsingular if and only if $(G_i,\bt_i)$ is nonsingular for $i=1,2$.

\ble\lb{m}
 Assume $(G,\bt)$ to be the orthogonal direct sum of the symplectic abelian groups $(G_i,\bt_i),\ i=1,2$. Let $\xi_i\in Z^2(G_i, F^\times)$ be a 2-cocycle corresponding to $\bt_i$ and $F^{\xi_i} G_i$ be the twisted group algebra of $(G_i,\bt_i)$ for $i=1,2$, then

(1) $\xi\in Z^2(G, F^\times)$ defined by $\xi((g_1,g_2),(h_1,h_2))=\xi_1(g_1,h_1) \xi_2(g_2,h_2)$ is a 2-cocycle corresponding to $\bt$.

(2) $F^\xi G\cong F^{\xi_1} G_1\ot F^{\xi_2} G_2$ is the twisted group algebra of $(G,\bt)$.
\ele

\bp (1) is directly verified.

(2) It is easy to check that the map $$F^{\xi_1} G_1\ot F^{\xi_2} G_2\rt F^\xi G, u_{g_1}\ot u_{g_2}\mapsto u_{(g_1,g_2)}   $$ is an algebra isomorphism.

\ep
Note that orthogonal direct sum can be defined for any finite number of symplectic abelian groups and above lemma can be generalized to that case similarly.

 \bpr\lb{34}
Let $(G,\bt)$ be a finite symplectic abelian group with $G=\op_{i=1}^k \mathbb{Z} a_i$ and $ord(a_i)=n_i$. Let $\xi\in L^2(G,F^\times)$ be defined by
\be\lb{341} \xi(a_i,a_j)=\begin{cases} 1 & \text{if\ $i\leq j$};\\
\bt(a_i,a_j) & \text{if\ $i> j$.}\end{cases} \ee Denote $u_{a_i}$ by $\overline{a_i}$. Then

(1) $\Psi(\xi)=\bt$.

(2) Let $A$  be the associative algebra over $F$ generated by $\{x_i|i=1,\cdots,k\}$ subject to the relations

$(a) x_i x_j=\bt(a_i,a_j) x_j x_i,i,j=1,\cdots,k $

and

$(b) x_i^{n_i}=1,i=1,\cdots,k$.

Then $\phi:A\rt F^\xi G, x_1^{l_1}\cdots x_k^{l_k}\mapsto \overline{a_1}^{l_1}\cdots \overline{a_k}^{l_k}$ defines an algebra isomorphism. In particular $ F^\xi G$ is generated by $\overline{a_i}, i=1,\cdots,k$, and \[\overline{a_1}^{l_1}\cdots \overline{a_k}^{l_k}=\overline{l_1 a_1+\cdots+l_k a_k}.\]
\epr
\bp
(1) is directly verified.

(2) Since $G$ is finite, it is easy to see that $F^\xi G$ is generated by $\overline{a_i}$. We first verify that $\overline{a_1}, \cdots, \overline{a_k}$ satisfies the relations (a) and (b).

  (a) is clear since $\overline{a_i}\ \overline{a_j}=\bt(a_i,a_j) \overline{a_j}\ \overline{a_i},i,j=1,\cdots,k $. For  (b), $\overline{a_i}^2=\xi(a_i,a_i)\overline{2a_i}=\overline{2a_i}$ and by induction one can show that  \[\overline{a_i}^k=\xi(a_i,(k-1)a_i)\overline{k a_i}=\overline{k a_i}\] for any positive integer $k$. In particular, let $k=n_i$ and one has $\overline{a_i}^{n_i}=\overline{0}=1$ since $ord(a_i)=n_i$.

Then the algebra homomorphism from the free algebra generated by $x_1,\cdots,x_k$ to $F^\xi G$ that maps $x_{i}$ to $\overline{a_i}$ descends to the algebra homomorphism $\phi:A\rt F^\xi G$ that maps $x_1^{l_1}\cdots x_k^{l_k}$ to $\overline{a_1}^{l_1}\cdots \overline{a_k}^{l_k}$.

%One has \bee  x_1^{l_1}\cdots x_k^{l_k}\cdot x_1^{j_1}\cdots x_k^{j_k}=\prod_{k\geq s>t\geq 1} \bt(a_s,a_t)^{l_s j_t}\  x_1^{l_1+j_1}\cdots x_k^{l_k+j_k}. \eee
%and
%\begin{align*}\overline{a_1}^{l_1} \cdots  \overline{a_k}^{l_k} \cdot \overline{a_1}^{j_1} \cdots  \overline{a_k}^{j_k}=&\prod_{k\geq s>t\geq 1} \bt(a_s,a_t)^{l_s j_t}\overline{a_1}^{l_1}\overline{a_1}^{j_1} \cdots  \overline{a_k}^{l_k}\overline{a_k}^{j_k}\\=&\prod_{k\geq s>t\geq 1} \bt(a_s,a_t)^{l_s j_t}\overline{a_1}^{l_1+j_1} \cdots  \overline{a_k}^{l_k+j_k} \end{align*}

Since $F^\xi G$ is generated by $\overline{a_i}$, $\phi$ is a surjective algebra homomorphism, which is also an isomorphism since by relation (a) $A$ is spanned by $x_1^{l_1}\cdots x_k^{l_k}, l_i=0,1,\cdots, n_i-1$. The last equation follows from
(\ref{341}).
\ep

Let $A$ be a (maybe infinite dimensional) associative algebra with a multiplicity free division $G$-grading $\Ga$, and identify $A_0=F\cdot 1$ with $F$.
\ble
(1) One has  $A_g A_h=A_{g+h}$ for all $g,h\in \Supp\ \Ga$.

(2) If $g\in \Supp\ \Ga$ and $A_g=F\cdot a_g$ then $A_{-g}=F\cdot a_g^{-1}$.

(3) $\Supp\ \Ga=G$.
\ele
\bp
(1) For nonzero $a_g\in A_g, a_h\in A_h$, $0\neq a_g a_h\in A_{g+h}$ as $a_g, a_h$ are both invertible. As $dim\ A_{g+h}=1$, $A_{g+h}=F\cdot  a_g a_h=A_g A_h$.

(2) As $A$ is graded, $a_g^{-1}$ can be uniquely written as $a_{-g}+b$, where $a_{-g}$ is the homogeneous component of $a_g^{-1}$ with degree $-g$ and $b=a_g^{-1}-a_{-g}$ is the sum of  the homogeneous components of $a_g^{-1}$ with degree $\neq -g$. $a_g^{-1}\cdot a_g=1$ implies that $a_{-g}\cdot a_g=1$ and $b\cdot a_g=0$. Since  $a_g$ is invertible,  $a_{-g}=a_g^{-1}$ and $b=0$. Thus $a_g^{-1}\in A_{-g}$ and $A_{-g}=F\cdot a_g^{-1}$as $dim\ A_{-g}=1$.

(3) By (1) and (2) we see that if  $g,h\in \Supp\ \Ga$ then $g+h, -g\in \Supp\ \Ga$. Then $\Supp\ \Ga$ is a subgroup of $G$. But as it is assumed that $\Supp\ \Ga$ generates $G$, one has $G=\Supp\ \Ga$.

\ep

\bpr\lb{3a}
Let $A$ be a (maybe infinite dimensional) associative algebra with a multiplicity free division $G$-grading $\Ga$. For any $g\in G$ let $a_g\in A_g$ be a non-zero element and define $\xi:G\times G\rt F^\times, \xi(g,h)=a_g a_h a_{g+h}^{-1}$. Then $\xi\in Z^2(G,F^\times)$ and $A$ is isomorphic to the twisted group algebra $F^\xi G$ as a $G$-graded algebra. The bicharacter $\bt=\Psi(\xi)\in\w^2(G,F^\times)$ is as follows:    \be \bt(g,h)=\xi(g,h)\xi(g,h)^{-1}=a_g a_h (a_g)^{-1} (a_h)^{-1}, \ for \ any\ g,h\in G.\ee It is clear that $\bt(g,h)$ does not depend on the choice of $a_g,a_h$.
\epr
The proof is direct and  we omit the details. The linear map \[A\rt F^\xi G, a_g\mapsto u_g\] is an isomorphism of $G$-graded associative algebras.

\section{Skew root systems and corresponding Lie algebras and Jordan algebras}
 \setcounter{equation}{0}\setcounter{theorem}{0}

Now we  define skew root systems of Lie type and of Jordan type.

 Let $(G,\bt)$ be a symplectic abelian group with $G$ a finitely generated additive abelian group.
\begin{defi}\lb{}
Assume $\bt$ to be nontrivial, i.e., $\Rad(\bt)\neq G$.
If a subset $R$ of $G$ satisfies
\begin{itemize}
  \item \textrm{(SRSL0)}. $R\subseteq G\setminus \Rad(\bt)$ and $R$ generates $G$.
 \item \textrm{(SRSL1)}. If $a\in R$, then $-a\in R$.
  \item \textrm{(SRSL2)}. If $\bt(a,b)\neq 1$, then $a+b\in R$.

\end{itemize} Then we say that $(G,\bt,R)$ is a \textbf{skew root system of Lie type}, or $R$ is a \textbf{skew root system of Lie type} in $(G,\bt)$.
  \end{defi}

\begin{defi}\lb{}
If a subset $R$ of $G$ satisfies
\begin{itemize}
  \item \textrm{(SRSJ0)}. $R\subseteq G$ and $R$ generates $G$.
 \item \textrm{(SRSJ1)}. If $a\in R$, then $-a\in R$.
  \item \textrm{(SRSJ2)}. If $\bt(a,b)\neq -1$, then $a+b\in R$.
\end{itemize}
Then we say that $(G,\bt,R)$ is a \textbf{skew root system of Jordan type}, or $R$ is a \textbf{skew root system of Jordan type} in $(G,\bt)$.
  \end{defi}

When we refer to $(G,\bt,R)$ as a skew root system, it maybe either of Lie type or of Jordan type.
 \ble
  If $(G,\bt,R)$ is a {skew root system of Jordan type}, then $0\in R$.
  \ele
 \bp
 For $a\in R$, $-a\in R$ by (SRSJ1). Then $\bt(a,-a)=1\neq -1$ since the characteristic of $F$ is 0, thus $0=a+(-a)\in R$ by (SRSJ2).
 \ep

   \begin{rem}
  If $(G,\bt,R)$  is a skew root system of Lie type, then $0\notin R$. If $(G,\bt,R)$  is a skew root system of Jordan type, then $0\in R$. So a skew root system cannot be of Lie type and of Jordan type at the same time.
  \end{rem}

  Two skew root systems $(G_i,\bt_i,R_i)$ ($i=1,2$) of the same type are \textit{isomorphic}, denoted $(G_1,\bt_1,R_1)\cong (G_2,\bt_2,R_2)$, if there is a group isomorphism $\varphi:G_1\rt G_2$ preserving the respective alternating bicharacters and $\varphi(R_1)=R_2$.

If $A$ is an associative algebra, then $A$ can be viewed as a Lie algebra with the Lie bracket $[a,b]=ab-ba$ and denote this Lie algebra by $A^-$. Similarly if $A$ is an associative algebra, then $A$ can be viewed as a Jordan algebra with the Jordan product $a\circ b=\frac{1}{2}(ab+ba)$ and denote this Jordan algebra by $A^+$. Because the group $G$ is abelian, one has
\ble\lb{k}
If an associative algebra $A$ has a $G$-grading, then this grading is  a  $G$-grading on the Lie algebra $A^-$, and is also a  $G$-grading on the Jordan algebra $A^+$.
\ele
\bp
For any $g,h\in G$ and $a\in A_g,b\in A_h$, one has $ab,ba\in A_{g+h}$. Thus  $[a,b],a\circ b\in A_{g+h}$. So $[A_g,A_h]\subseteq A_{g+h}$ and $A_g\circ A_h\subseteq  A_{g+h}$ and the lemma follows.
\ep

Given a symplectic abelian group $(G,\bt)$, fix some $\xi\in Z^2(G,F^\times)$ in the preimage of $\bt$ under $\Psi: Z^2(G,F^\times)\rt \w^2(G,F^\times)$. Then one associates $(G,\bt)$ with a twisted group algebra $F^\xi G$ with the multiplicity free division grading \[F^\xi G=\oplus_{a\in G}F u_a.\]

\ble\lb{461}
(1)Given a skew root system $(G,\bt,R)$ of Lie type, there corresponds to it a $G$-graded Lie algebra \be\lb{41} L(R)=\oplus_{a\in R}F u_a,\ee which is a $G$-graded Lie subalgebra of $(F^\xi G)^-$ with \be\lb{46}[u_a,u_b]=u_a u_b-u_b u_a=(\xi(a,b)-\xi(b,a))u_{a+b}.\ee

(2)Similarly given a skew root system $(G,\bt,R)$ of Jordan type, there corresponds to it a $G$-graded Jordan algebra \be\lb{42} J(R)=\oplus_{a\in R}F u_a,\ee which is a $G$-graded Jordan subalgebra of $(F^\xi G)^+$ with \be\lb{47}u_a \circ u_b=\frac{1}{2}(u_a u_b+u_b u_a)=\frac{1}{2}(\xi(a,b)+\xi(b,a))u_{a+b}.\ee
\ele
\bp (1) $L(R)$ is closed under the Lie bracket (\ref{46}) is equivalent to the condition if $a,b\in R$ and $\xi(a,b)-\xi(b,a)\neq 0$ then $a+b\in R$. But $\xi(a,b)-\xi(b,a)\neq 0$ if and only if $\bt(a,b)\neq 1$. So the axiom (SRSL2) guarantees that $L(R)$ is a $G$-graded Lie subalgebra of $(F^\xi G)^-$.

(2) $J(R)$ is closed under the Jordan multiplication (\ref{47}) is equivalent to the condition if $a,b\in R$ and $\xi(a,b)+\xi(b,a)\neq 0$ then $a+b\in R$. But $\xi(a,b)+\xi(b,a)\neq 0$ if and only if $\bt(a,b)\neq -1$. So the axiom (SRSJ2) guarantees that $J(R)$ is a $G$-graded Jordan subalgebra of $(F^\xi G)^+$.
 \ep
Let $(G,\bt,R)$ be a skew root system . It is called \textit{reduced} if the alternating bicharacter $\bt$ is nonsingular. Next lemma implies that any skew root system  corresponds to a reduced skew root system .
\ble\lb{47}
Assume that $(G,\bt,R)$ is a skew root system and $H$ is any subgroup of $\Rad(\bt)$. Let $\overline{G}=G/H$, $\overline{R}$ be the image of $R$ in $\overline{G}$, and $\overline{\bt}$ be the alternating bicharacter on $\overline{G}$ induced by $\bt$, i.e., $\overline{\bt}(\overline{g},\overline{h})=\bt(g,h)$ for any $\overline{g},\overline{h}\in \overline{G}$.
Then

(1) $(\overline{G},\overline{\bt},\overline{R})$ is a skew root system (of the same type as that of $(G,\bt,R)$). It is reduced if $H=\Rad(\bt)$ and is referred as the corresponding reduced skew root system of $(G,\bt,R)$.

(2) There is a surjective  algebra homomorphism from $L(R)$ to $L(\overline{R})$ that maps $L(R)_g$ to $L(\overline{R})_{\overline{g}}$ for any $g\in R$ if $(G,\bt,R)$ is of Lie type.  Similarly,  There is a surjective  algebra homomorphism from $J(R)$ to $J(\overline{R})$ that maps $J(R)_g$ to $J(\overline{R})_{\overline{g}}$ for any $g\in R$ if $(G,\bt,R)$ is of Jordan type.
%If $L(R)$ is simple then this is an isomorphism of Lie algebras, and in this case the homogeneous components of $L(R)$ and $L(\overline{R})$ are the same.
\ele
\bp
(1) It is directly verified that $(\overline{G},\overline{\bt},\overline{R})$ satisfies the axioms of skew root systems thus is also a skew root system . Note that $\Rad(\overline{\bt})=\Rad(\bt)/H$. So if $H=\Rad(\bt)$ then it is a reduced skew root system.

(2) We will prove it in the case $(G,\bt,R)$ is of Lie type, and the proof in the Jordan case is similar.

 Let $\overline{\xi}\in Hom_\mathbb{Z}(\overline{G}\ot \overline{G}, F^\times)$ such that $\Psi(\overline{\xi})=\overline{\bt}$. Let  ${\xi}\in Hom_\mathbb{Z}({G}\ot {G}, F^\times)$ be the pull-back of $\overline{\xi}$, i.e., for any $g,h\in G$, $\xi(g,h)=\overline{\xi}(\overline{g},\overline{h})$. Then it is directly verified that $\Psi({\xi})={\bt}$.

Then $\phi:L(R)\rt L(\overline{R}), u_g\mapsto u_{\overline{g}}$ defines a surjective Lie algebra homomorphism, as $$[u_g,u_h]=(\xi(g,h)-\xi(h,g))u_{g+h}$$ and $$[u_{\overline{g}},u_{\overline{h}}]=(\overline{\xi}(\overline{g},\overline{h})-\overline{\xi}(\overline{h},\overline{g}))u_{\overline{g}+\overline{h}}
=(\xi(g,h)-\xi(h,g))u_{\overline{g+h}}.$$

\ep

A skew root system $(G,\bt,R)$ of Lie type is called \textit{decomposable} if $G$ is an orthogonal sum of two proper subgroups $G_1$ and $G_2$, $R$ is a union of two  orthogonal subsets $R_1$ and $R_2$ (i.e., $\bt(a_i,a_2)=1$ for any $a_i\in R_i$ with $i=1,2$)  , and $R_i$ is a skew root system in $G_i$ for $i=1,2$. We say that $(G,\bt,R)$ is the \textit{direct sum} of $(G_i,\bt_i,R_i)$ for $i=1,2$. It is clear that in this case $L(R)=L(R_1)\op L(R_2)$ is a direct sum of ideals $L(R_1)$ and $L(R_2)$. Similarly a skew root system of Lie type could be a direct sum of finitely many skew root systems of Lie type. The skew root system $(G,\bt,R)$ of Lie type is \textit{indecomposable} if it is not decomposable.

For a skew root system  $(G,\bt,R)$, let $\epsilon=1$ if it is of Lie type and $\epsilon=-1$ if it is of Jordan type. Define a simple graph $(R, E)$ on the set $R$ (where $E$ is set of edges with vertices in $R$) as follows: For any $a,b\in R$ with $a\neq b$, $\{a,b\}\in E$ if and only if $\bt(a,b)\neq \epsilon$. We refer to $(R, E)$ as the \textit{graph} of $(G,\bt,R)$. Note that (SRSL2) and (SRSJ2) implies that if $\{a,b\}\in E$ then $a+b\in R$. If $(G,\bt,R)$ is of Jordan type, then $\{0,a\}\in E$ for any $a\in R\setminus \{0\}$ thus $(R, E)$ is connected. But if $(G,\bt,R)$ is of Lie type, then $(R, E)$ may not be connected.

\ble
Let $(G,\bt,R)$ be a reduced skew root system  of Lie type and $(R, E)$ be the corresponding graph. Let $R_1$ be a connected component of $(R, E)$. Then $|R_1|>1$ and $R_1=-R_1$.
\ele
\bp
If $|R_1|=1$ and assume $R_1=\{a\}$, then $a\perp a$ and $a\perp (R \setminus R_1)$. So $a\perp R$ and $a\in \Rad(\bt)$ as $G$ is generated by $R$. But $(G,\bt)$ is nonsingular, thus $a=0$ which contradicts to $(SRSL0)$. So $|R_1|>1$.

      For any $a\in R_1$, there is some $b\in R_1$ with $\bt(a,b)\neq 1$. Then $\bt(-a,b)=\bt(a,b)^{-1}\neq 1$. So $-a\in R_1$ thus $R_1=-R_1$.
\ep
\ble
Let $(G,\bt,R)$ be a reduced skew root system  of Lie type and $(R, E)$ be the corresponding graph. Assume that the graph is disconnected. Let $R_1$ be a connected component of $(R, E)$ and $R_2=R\setminus R_1$. For $i=1,2$ Let $G_i$ be the subgroup of $G$ generated by $R_i$ and $\bt_i=\bt|_{G_i}$ be the alternating bicharacter on $G_i$ induced from $\bt$. Then

(1) $(G,\bt)$ is the orthogonal direct sum of $(G_i,\bt_i)$ for $i=1,2$, i.e., $G_1\cap G_2=0$, $G_1+G_2=G$ and $G_1\perp G_2$. Moreover, $(G_i,\bt_i)$ is also nonsingular.

(2) One has $R_i$ is a skew root system  of Lie type in $(G_i,\bt_i)$ for $i=1,2$. The skew root system  $(G,\bt,R)$ is the {direct sum} of $(G_i,\bt_i,R_i)$ for $i=1,2$ and is decomposable.
\ele
\bp
(1)  Assume $g\in G_1\cap G_2$. Then $g\perp R_2$ since $g$ is a finite sum of elements in $R_1$ and $R_1\perp R_2$. Symmetrically $g\perp R_1$. Then $g\perp R$ and $g\in \Rad(\bt)$. So $g=0$ and $G_1\cap G_2=0$. As $G_1+G_2$ is a subgroup of $G$ containing $R=R_1\cup R_2$, and $R$ generates $G$, one has $G_1+G_2=G$. The fact $R_1\perp R_2$ implies that $G_1\perp G_2$. So $(G,\bt)$ is the orthogonal direct sum of $(G_i,\bt_i)$ for $i=1,2$. The last assertion is clear.

(2) Assume  $i=1$ or 2. It is clear that $(G_i,\bt_i,R_i)$ satisfies (SRSL0), and also satisfies (SRSL1) by last lemma. Next we verify that $(G_i,\bt_i,R_i)$ satisfies (SRSL2). If $a,b\in R_i$ and $\bt(a,b)\neq 1$, then $a+b\in R$. As $\bt(a,a+b)=\bt(a,b)\neq 1$, $a+b$ is in the same connected component of $a$ and $b$ thus $a+b\in R_i$. It is clear that $(G,\bt,R)$ is the {direct sum} of $(G_i,\bt_i,R_i)$ for $i=1,2$ and is decomposable.

\ep

\bco\lb{49}
Let $(G,\bt,R)$ be a reduced skew root system  of Lie type and $(R, E)$ be the corresponding graph. Then

(1) $(G,\bt,R)$ is indecomposable if and only if $(R, E)$ is connected.

(2) Assume that $(R, E)$ has finite many connected components, say $R_i (i=1,\cdots, n)$. Let $G_i$ be the subgroup of $G$ generated by $R_i$. Then $(G,\bt,R)$ is the direct sum of the reduced and indecomposable skew root system  $(G_i,\bt_i,R_i)$, where $\bt_i=\bt|_{G_i}$, and $L(R)=\op_i L(R_i)$.

\eco
\bp
(1) Last lemma showed if $(R, E)$ is disconnected then $(G,\bt,R)$ is decomposable. Conversely, assume $(G,\bt,R)$ to be decomposable. By definition, $G$ is an orthogonal sum of two proper subgroups $G_1$ and $G_2$, $R$ is a union of two  orthogonal subsets $R_1$ and $R_2$, and $R_i$ is a skew root system in $G_i$ for $i=1,2$. Then there is no edge between $R_1$ and $R_2$ thus $(R, E)$ is disconnected.

  (2) This follows from last lemma.

\ep
In view of Corollary \ref{49} (1), any skew root system of Jordan type is called \textit{indecomposable} as its graph is connected.
\ble
 Assume that $(G,\bt,R)$ is a skew root system  of Lie type that is reduced and indecomposable.

  (1) Let $I$ be an ideal of $L(R)$ containing $u_{a}$ for some $a\in R$. Then $I=L(R)$.

  (2) One has that $L(R)$ is graded-simple.
\ele
\bp
(1) If $R=\{a\}$ then $G$ (generated by $R$) is cyclic which is impossible as in this case $\Rad(\bt)=G$. Since $R$ is indecomposable, the graph of $R$ is connected. For any $b\in R$, there is a path $a=a_0,a_1,\cdots,a_n=b$ linking $a$ and $b$, i.e., $\{a_i,a_{i+1}\}\in E$. We will show that $u_{a_i}\in I$ implies $u_{a_{i+1}}\in I$, from which we conclude that $I=L(R)$.

       One has $a_{i+1}-a_i\in R$ since $\bt(a_{i+1},-a_i)\neq 1$. Choose $\xi\in L^2(G,F^\times)$ corresponding to $\bt$. Then \be \begin{split} [u_{a_i},u_{a_{i+1}-a_i}]&=(\xi(a_i,{a_{i+1}-a_i})-\xi({a_{i+1}-a_i},a_i))u_{a_{i+1}}\\
       &=\xi(a_i,a_i)^{-1}(\xi(a_i,{a_{i+1}})-\xi({a_{i+1}},a_i))u_{a_{i+1}}\\
       &=\xi(a_i,a_i)^{-1}\xi({a_{i+1}},a_i)(\bt(a_i,{a_{i+1}})-1)u_{a_{i+1}}\end{split}\ee is a nonzero multiple of $u_{a_{i+1}}$ as $\bt(a_i,{a_{i+1}})\neq 1 $.

(2) follows from (1).
\ep
 \bthm\lb{411}Assume $(G,\bt,R)$ is a reduced skew root system  of Lie type.

   (1)One has that $R$ is indecomposable if and only if $L(R)$ is simple.

   (2)If $(G,\bt,R)$ is a direct sum of indecomposable skew root system  $(G_i,\bt_i,R_i)$ for $i=1,\cdots,n$, then $L(R)$ is the direct sum of simple Lie algebras $L(R_i)$.
  \ethm
 \bp
 We will show that if $R$ is indecomposable then $L(R)$ is simple. Then (2) will follows from this fact and Corollary \ref{49} (2). Next by (2) we get that if $L(R)$ is simple then $R$ is indecomposable, which complete the proof  of (1).

 Define the support and length for any $0\neq z\in L(R)$ as follows.
The element $z$ can be uniquely written as a finite sum $\sum_{\al\in S} \lambda_\al u_\al$ where $S\subseteq R$ and $0\neq\lambda_\al\in F$ for any $\al\in S$. Define the support $supp(z)$ of $z$ to be $S$, and the length $l(z)$ of $z$ to be the cardinality of $S$.

Let $I$ be a nonzero ideal of $L(R)$. Let $x$ be a nonzero element in $I$ with minimal length.

Case 1. $l(x)=1$. Then last lemma implies that $I=L(R)$.

Case 2. $l(x)>1$. Let $S$ be its support. There are 2 subcases.

Case 2.1. There exist $a,b\in S$ with $\bt(a,b)\neq 1$. As $[u_a, u_a]=0, [u_a, u_b]\neq 0$,  $y=[u_a, x]\in I$ is nonzero and $l(y)<l(x)$, which contradicts to that $l(x)$ is minimal.

Case 2.2. For any $a,b\in S$, $\bt(a,b)=1$. Choose $a,b\in S$. As $(G,\bt)$ is nonsingular and $R$ generates $G$, there exists some $c\in R$ with $\bt(a-b,c)\neq 1$. As $\bt(a-b,c)=\bt(a,c)/\bt(b,c)\neq 1$, $\bt(a,c)\neq\bt(b,c)$. There are still 2 subcases.

Case 2.2.1. One of $\bt(a,c)$ and  $\bt(b,c)$ is 1. Then  $y=[u_c, x]\in I$ is nonzero and $l(y)<l(x)$, which contradicts to that $l(x)$ is minimal.

Case 2.2.2. $\bt(a,c)\neq 1$ and  $\bt(b,c)\neq 1$. Let $y=[u_c, x]\in I$. Then $\{a+c,b+c\}\subseteq supp(y)\subseteq\{d+c|d\in\ S\}$ and $supp(y)\leq supp(x)$. As $\bt(a+c,b+c)=\bt(a,c)\bt(c,b)=\bt(a-b,c)\neq 1$, if we replace $x$ by $y$ we get to Case 2.1, which will also lead to a contradiction.

So one always has $I=L(R)$ and $L(R)$ is simple.
\ep
\bpr\lb{412}
 (1)  If $(G,\bt,R)$ is  skew root system  of Jordan type, then $J(R)$ is graded-simple.

 (2) If $(G,\bt,R)$ is a reduced skew root system  of Jordan type, and $G$ is an elementary abelian 2-group, then $J(R)$ is simple.
 \epr
 \bp
    (1)   Assume $I$ to be a non-zero graded ideal of $J(R)$. Then $u_{g}\in I$ for some $g\in R$.  Choose $\xi\in L^2(G,F^\times)$ corresponding to $\bt$. One has
   \be \lb{i}\begin{split} u_{g}\circ u_{-g}&=\frac{1}{2}(\xi(g,-g)+\xi(-g,g))u_{0}\\&=\frac{1}{2}\xi(-g,g)(\bt(g,-g)+1)u_{0}\\&=\xi(-g,g)u_{0}.\end{split}\ee

    So $u_{0}=\xi(g,g)u_{g}\circ u_{-g}\in I$. Since $u_{0}$ is the identity of $J(R)$, $I=J(R)$.

    (2) Define the support and length for any nonzero element in $J(R)$ as in last theorem. Let $I$ be a non-zero ideal of $J(R)$ and $x$ be a nonzero element in $I$ with minimal length.

Case 1. $l(x)=1$. Then the proof of (1) in fact implies that $I=L(R)$.

Case 2. $l(x)>1$. Let $S$ be its support. There are 2 subcases.

Case 2.1. There exist $a,b\in S$ with $\bt(a,b)=-1$. As $u_a\circ u_a\neq 0, u_a\circ u_b=0$,  $y=u_a\circ x\in I$ is nonzero and $l(y)<l(x)$, which contradicts to that $l(x)$ is minimal.

Case 2.2. For any $a,b\in S$, $\bt(a,b)\neq -1$. Choose $a,b\in S$. As $(G,\bt)$ is nonsingular and $R$ generates $G$, there exists some $c\in R$ with $\bt(a-b,c)\neq 1$. Then $\bt(a,c)\neq\bt(b,c)$. As $G$ is an elementary abelian 2-group, $\bt(g,h)=\pm 1$ for any $g,h\in G$. So one of $\bt(a,c)$ and $\bt(b,c)$ is 1 and the other  is $-1$. Then $y=u_c\circ x\in I$ is nonzero and $l(y)<l(x)$, which contradicts to that $l(x)$ is minimal.

    So in the case $G$ is an elementary abelian 2-group one always has $I=J(R)$ thus $J(R)$ is simple.

\ep

 Now we assume that $G$ is a finite group.
  \begin{prop}\lb{413} If $G$ is a finite group and $(G,\bt,R)$ is a skew root system of Lie type, then $L(R)$ is a semisimple Lie algebra. The grading on $L(R)$ is a multiplicity free semisimple grading.
\end{prop}

 \bp
  Choose $\xi\in L^2(G,F^\times)$ with $\Psi(\xi)=\bt$, i.e., $\bt(a,b)=\xi(a,b)\xi(b,a)^{-1}.$ Then the $G$-graded algebra $L(R)$ is constructed as in (\ref{41}).

  Let $(,)$ be the Killing form on $L(R)$. It is clear that \be\lb{4a} (u_a,u_b)=0, if\ a+b\neq 0\ee as $ad_{u_{a}}\cdot ad_{u_b}$ is nilpotent.
 As $G$ is finite, $\xi(a,b)$ and $\bt(a,b)$ is always a root of unity.
 One has  \be\lb{o}\begin{split}(u_{-a},u_{a})&=tr(ad_{u_{-a}}\cdot ad_{u_a} )\\&=\sum_{b\in R} (\xi(a,b)-\xi(b,a))(\xi(-a,a+b)-\xi(a+b,-a))\\&=\sum_{b\in R}\xi(a,-a) (2-\xi(a,b)\xi(b,-a)-\xi(b,a)\xi(-a,b))\\&=\xi(a,-a)\sum_{b\in R}(2-\bt(a,b)-\bt(b,a))\\&=\xi(a,-a)\sum_{b\in R}(2-(\bt(a,b)+\bt(a,b)^{-1})).\end{split}\ee  As $F$ is algebraically closed with characteristic 0, $F$ contains the algebraic closure of the field $Q$ of rational numbers and
\[2-(\bt(a,b)+\bt(a,b)^{-1}\geq 0,\] with equality holds if and only if  $\bt(a,b)=1$.

For any $a\in R\subseteq G\setminus \Rad(\bt)$, there is some $b\in R$ with $\bt(a,b)\neq 1$. So  \be\lb{4b}(u_a,u_{-a})\neq 0,\ for \ all\ a\in R.\ee Then (\ref{4a}) and (\ref{4b}) implies that $(,)$ is nonsingular and $L(R)$ is semisimple. It is clear that each $u_a$ is semismiple, thus the grading on $L(R)$ is a multiplicity free semisimple grading.
 \ep

 \bpr\lb{414}
 If $G$ is a finite group and $(G,\bt,R)$ is a skew root system of Jordan type, then $J(R)$ is a semisimple $G$-graded Jordan algebra. The grading on $J(R)$ is a multiplicity free division grading.
 \epr
 \bp Choose $\xi\in L^2(G,F^\times)$ corresponding to $\bt$. For any $u\in J(R)$, let $L(u)$ be the linear transformation on $J(R)$ defined by left multiplication by $u$.
 Let $\tau$ be the canonical symmetric bilinear form on $J(R)$, i.e., for any $u,v\in J(R)$, $$\tau(u,v)=tr\ L(u\circ v).$$ One knows that $\tau$ is associative. Next we will show that $\tau$ is nonsingular which implies that $J(R)$ is semisimple.

 Note that if $g,h\in G$ and $g+h\neq 0$ then $L(u_{g}\circ u_{h})$ is nilpotent and \be\lb{e}\tau(u_{g},u_{h})=0.\ee By (\ref{i}) one has
\be\lb{h} \begin{split} \tau(u_{g},u_{-g})&=tr\ L(u_{g}\circ  u_{-g})\\&=tr\ L(\xi(-g,g)u_{0})\\&=\xi(-g,g) dim\ J\neq 0 \end{split}\ee since the characteristic of $F$ is 0.  By (\ref{e}) and (\ref{h}) it is easy to deduce that $\tau$ is nonsingular.

It is easy to see that $u_{g}$ is invertible in the Jordan algebra $J(R)$ and its inverse is $\xi(g,g)u_{-g}$ by (\ref{i}). Thus the $G$-grading on $J(R)$ is a multiplicity free division grading.
  \ep

 \section{Constructions of skew root system  and corresponding algebras}
 \setcounter{equation}{0}\setcounter{theorem}{0}

The task to classify skew root systems  $(G,\bt,R)$ for finitely generated abelian group $G$  seems to be hard even if we assume that the skew root system  is reduced and indecomposable. In this section we will construct 3 families of  skew root systems (which may be nonreduced) of Lie type and of Jordan type respectively in finite  symplectic abelian groups, and identify the corresponding Lie and Jondan algebras respectively, most of which are simple.

In the following 3 cases respectively, we will describe the symplectic abelian group  $(G,\bt)$ and the corresponding twisted group algebra first. Then we construct the skew root system  $R^-$ of Lie type and the skew root system  $R^+$ of Jordan type and then identify $L(R^-)$  and $L(R^+)$ respectively. \bigskip

\textbf{5.1.  $(G,\bt)$ is a finite nonsingular symplectic abelian group with $R^-=G\setminus \{0\}$ and $R^+=G$.}\bigskip

 Nonsingular symplectic abelian groups $(G,\bt)$ are classified and the corresponding twisted group algebras $F^\xi(G)$ are simple, see Theorem 1.8 and Theorem 2.12 of \cite{ka3}. In this case it is clear that $R^-=G\setminus \{0\}$ satisfies the axioms of a skew root system  of Lie type and $R^+=G$ satisfies the axioms of a skew root system  of Jordan type.

For any $n\in  \mathbb{Z}_+$, let $ \mathbb{Z}_n= \mathbb{Z}/n\mathbb{Z}=\{\bar{0},\bar{1},\cdots,\overline{n-1}\}$. For simplicity
 we will just use an integer $i$ to denote $\bar{i}$.

First let $G=\mathbb{Z}_n\times \mathbb{Z}_n$ for some $n>1$. Let $\varepsilon\in F$ be a primitive $n$-th root of unity. For any $(i,j),(s,t)\in \mathbb{Z}_n\times \mathbb{Z}_n$, let $\bt((i,j),(s,t))=\varepsilon^{it-js}$ (which is defined by $\bt(e_1,e_2)=\varepsilon$ with $e_1=(1,0)$ and $e_2=(0,1)$). It is a nonsingular alternating bicharacter on $G$. Any  finite nonsingular symplectic abelian group is a direct sum of such $(G,\bt)$'s (where the order of $G$'s may be different).

Let $\xi\in L^2(G, F^\times)$ be defined by $\xi((i,j),(s,t))=\varepsilon^{-js}$, which is defined as in (\ref{341}) of Proposition \ref{34}. Then $\Psi(\xi)=\bt$.

One knows that $F^\xi(G)\cong M_n(F)$. A concrete isomorphism is as follows. Let $X_1, X_2\in M_n(F)$ be defined by
 \be \lb{y}X_1=\begin{pmatrix}
\varepsilon^{n-1} & 0 & \cdots & 0 \\
0 & \varepsilon^{n-2} & \cdots & 0 \\
 &  & \cdots &  \\
0 & 0 & \cdots & 1
\end{pmatrix},\quad X_2=\begin{pmatrix}
0 & 1 & 0 & \cdots & 0 \\
0 & 0 & 1 & \cdots & 0 \\
 &  &  & \cdots &  \\
0 & 0 & 0 & \cdots & 1 \\
1 & 0 & 0 & \cdots & 0
\end{pmatrix}. \ee
Then \be \lb{50}X_1 X_2=\varepsilon X_2 X_1,\ X_1^n=X_2^n=I_n\ee where $I_n$ is the $n\times n$ identity matrix.

One has \be\Ga: M_n(F)=\op_{(i,j)\in G} F\cdot X_1^i X_2^{j}\ee is a $G$-grading on $M_n(F)$, which is called the $\varepsilon$-grading on $M_n(F)$ in \cite{bsz}. In particular the associative algebra $M_n(F)$ is generated by $X_1$ and $X_2$ satisfying (\ref{50}), thus by proposition \ref{34}, \be \lb{54}\phi: F^\xi G\rt M_n(F), \overline{(i,j)}\mapsto X_1^i X_2^j\ee is a graded algebra isomorphism, where $\overline{(i,j)}=u_{(i,j)}$.

It is easy to see that the isomorphism $\phi$ induces a Lie algebra isomorphism \[\phi^-: L(R^-)\rt sl_n(F)\] and a Jordan algebra isomorphism \[\phi^+: J(R^+)\rt M_n(F)^+.\] %Note that $sl_n(F)$ and $M(n,F)^+$ is always simple if $n\geq 2$.

Now we deal with the general case. Assume $n=n_1\cdots n_k$ with each $n_i>1$ and \be \lb{} n_i|n_{i+1}, i=1,\cdots, k-1.\ee For $t=1,\cdots,k$, let $G_t=\mathbb{Z}_{n_t}\times \mathbb{Z}_{n_t}$ and $\varepsilon_t\in F$ be a primitive $n_t$-th root of unity. Let $\bt_t:G_t\times G_t\rt F^\times, \bt_t((i,j),(s,l))=\varepsilon_t^{il-js}$ be the alternating bicharacter on $G_t$ and $\xi_t:G_t\times G_t\rt F^\times, \xi_t((i,j),(s,l))=\varepsilon_t^{-js}$ as in the first case. Let $(G,\bt)$ be the  orthogonal direct sum of symplectic abelian groups $(G_1,\bt_1),\cdots, (G_k,\bt_k)$. Then $\bt$ is a nonsingular alternating bicharacter on $G$.

Let $\xi\in L^2(G,F^\times)$ be the product of $\xi_i$ defined in Lemma \ref{m} (1), then $\Psi(\xi)=\bt$. By Lemma \ref{m} $F^\xi G\cong F^{\xi_1} G_1\ot\cdots\ot F^{\xi_k} G_k\cong M(n_1,F)\ot\cdots\ot M(n_k,F)\cong M(n,F)$. By Lemma \ref{461}, the isomorphism induces  a Lie algebra isomorphism \[\phi^-: L(R^-)\rt sl_n(F)\] and a Jordan algebra isomorphism \[\phi^+: J(R^+)\rt M_n(F)^+.\] Note that $sl_n(F)$ and $M(n,F)^+$ is always simple if $n\geq 2$. By Theorem \ref{411} (1), the simpleness of $L(R^-)$ implies that $R^-$ is indecomposable.

 The main results in the subsection is summarized as follows.
\bpr
Let $(G,\bt)$ be a finite nonsingular symplectic abelian group. Then

(1) $R^-=G\setminus \{0\}$ is a reduced and indecomposable skew root system of Lie type in $(G,\bt)$ and $L(R^-)\cong sl_n(F)$ with $|G|=n^2$.

(2) $R^+=G$ is a reduced skew root system of Jordan type in $(G,\bt)$ and $J(R^+)\cong M_n(F)^+$ with $|G|=n^2$.

\epr

\textbf{5.2.  Skew root systems related to Clifford algebras}\bigskip

 Assume $n\geq 2$ and let $G=\mathbb{Z}_2^n$. Let $\{e_1,\cdots,e_n\}$ be the standard basis of $\mathbb{Z}_2^n$ and $\bt$ be the alternating bicharacter on $G$ defined by
\be\lb{551} \bt(e_i,e_j)=\begin{cases} 1, & \text{if\ i=j};\\
-1, & \text{if\ $i\neq j$.}\end{cases} \ee
The explicit formula  is \be\lb{58}\bt:G\times G\rt F^\times, \bt(a,b)=(-1)^{(\sum_{i=1}^{n} a_i)(\sum_{i=1}^{n} b_i)-\sum_{i=1}^{n} a_i b_i}.\ee Then

\bee \Rad(\bt)=\begin{cases} $\{0\}$, & \text{ n\ even};\\
$\{0,(1,... ,1)\}$,  & \text{ n\ odd}.\end{cases} \eee
Let
\bee \xi(e_i,e_j)=\begin{cases} 1, & \text{if\ $i\leq j$};\\
-1, & \text{if\ $i> j$.}\end{cases} \eee be the 2-cocycle on $G$ defined as in (\ref{341}) of Proposition \ref{34}. Then $\Psi(\xi)=\bt$.

The explicit formula for $\xi$ is \[\xi:G\times G\rt F^\times, \xi(a,b)=(-1)^{\sum_{1\leq j<i\leq n}(a_i b_j)}.\]

Let $V$ be an $n$-dimensional $F$-vector space with some nonsingular symmetric bilinear form and $\{v_1,\cdots,v_n\}$ be an orthonormal basis of $V$. Let $C(V)$ be the corresponding Clifford algebra. Then $C(V)$ is generated by $\{v_1,\cdots,v_n\}$ with relations
\bee v_i v_j=- v_j v_i, i\neq j; v_i^2=1, i=1,\cdots, n.\eee Then by Proposition \ref{34} one has the following isomorphism.
\ble\lb{4d}
The linear map $\phi: F^\xi G\rt C(V),  \overline{l_1 e_1+\cdots+l_n e_n}\mapsto v_1^{l_1}\cdots v_n^{l_n}$ defines an algebra isomorphism.
\ele

So $C(V)$ has the grading \be\lb{4c} C(V)=\bigoplus_{(l_1,\cdots,l_n)\in \mathbb{Z}_2^n} F\cdot v_1^{l_1}\cdots v_n^{l_n}.\ee

Let $R^+=\{0, e_i|1\leq i\leq n\}$.
Then it is easy to verify that $(G,\bt,R^+)$ is skew root system  of Jordan type.

$J(R^+)=F\cdot 1\op V$ is a $G$-graded Jordan subalgebra of $C(V)^+$, which is a simple Jordan algebra if $dim\ V>1$, i.e., $n\geq 2$, and is referred as the simple Jordan algebra of a nonsingular symmetric form. The $G$-grading on $J(R^+)$ is \be J(R^+)=F\cdot 1\op F\cdot v_1\op\cdots\op F\cdot v_n.\ee
%One knows that $Aut(J(R^+))$ is canonically identified with the orthogonal group $O(V)$. Then $\wh{G}\st Aut(J(R^+))$ is the group of diagonal matrices in $O(V)$ with diagonal entries $\pm 1$.

Let $R^-=\{e_i+e_j|1\leq i<j\leq n\}$ and $G_1\cong \mathbb{Z}_2^{n-1}$ be the subgroup of $G$ generated by $R$. Then it is easy to check that
\ble $\bt(e_i+e_j,e_s+e_t)=-1$ if and only if $|\{i,j\}\cap \{s,t\}|=1$.\ele

 Assume $n\geq 3$. Then this lemma implies that $(G_1,\bt,R^-)$ satisfies (SRSL2), since $(G_1,\bt,R^-)$ satisfies (SRSL0) and (SRSL1) obviously, it is a skew root system of Lie type for $n\geq 3$. It also follows from this lemma that the graph of $R^-$ is connected thus it is indecomposable.

 One has

\bee \Rad(\bt|_{G_1})=\begin{cases} $\{0\}$, & \text{ n\ odd};\\
$\{0,(1,... ,1)\}$,  & \text{ n\ even}.\end{cases} \eee
So $R^-$ is reduced if and only if $n$ is odd.

It is well known that graded subspace $L=\oplus_{1\leq i<j\leq n} F\cdot v_i v_j$ of $C(V)$ in (\ref{4c}) is a Lie algebra isomorphic to $so_n(F)$ under the usual Lie bracket. The isomorphism is
\bee L\rt so_n(F),  v_i v_j\mapsto 2(E_{ij}-E_{ij}). \eee Then the isomorphism in Lemma \ref{4d} induces an isomorphism of Lie algebras $L(R^-)$ and $L\cong so_n(F),$ which is simple if and only if $n\geq 3$ and $n\neq 4$.

 The main results in the subsection is summarized as follows.
\bpr
Let $(G,\bt)$ be the symplectic abelian group defined in (\ref{58}) with $G=\mathbb{Z}_2^n$ and $n\geq 2$. Then

(1)$R^+=\{0, e_i|1\leq i\leq n\}$ is a skew root system of Jordan type in $(G,\bt)$ which is reduced if and only if $n$ is even.
One has that $J(R^+)$ is isomorphic to the simple Jordan algebra of a nonsingular symmetric form.

(2)Let $R^-=\{e_i+e_j|1\leq i<j\leq n\}$ and $G_1\cong \mathbb{Z}_2^{n-1}$ be the subgroup of $G$ generated by $R^-$. Assume $n\geq 3$. Then  $(G_1,\bt,R^-)$ is a skew root system  of Lie type. It is always indecomposable, and is reduced if and only if $n$ is odd. One has that  $L(R^-)$ is isomorphic to the Lie algebra $so(n,F)$, which is simple if and only if $n\geq 3$ and $n\neq 4$.
\epr
This case can be generalized as follows. Let $G$ be the product of countably infinitely many  $\mathbb{Z}_2$ and the alternating bicharacter on $G$ is defined as in (\ref{551}). Then $R^-=\{e_i+e_j|1\leq i<j\}$ is a skew root system  of Lie type, and the Lie algebra $L(R^-)$ is an 'infinite dimensional special orthogonal Lie algebra'.

\bigskip

\textbf{ 5.3. Skew root systems associated with nonsingular quadratic forms over $\mathbb{F}_2$.}\bigskip

In the paper $\mathbb{F}_2$ denotes the field of two elements.

Let $q:\mathbb{F}_2^{n}\rt \mathbb{F}_2$ be a  quadratic form and \[\al_1:\mathbb{F}_2^{n}\times \mathbb{F}_2^{n}\rt \mathbb{F}_2, \al_1(a,b)=q(a+b)-q(a)-q(b)\] be its \textit{polarization}. Assume that $q$ is nonsingular, i.e.  $\Rad(q)=0$  where $\Rad(q)=\{a\in \Rad(\al_1)|q(a)=0.\}$.  One knows that up to isomorphism, there are two nonsingular quadratic forms on $\mathbb{F}_2^{n}$ for $n$ even, and there is one  nonsingular quadratic forms on $\mathbb{F}_2^{n}$ for $n$ odd. We will always identify the group $G=\mathbb{Z}_2^{n}$ with the additive group of $\mathbb{F}_2^{n}$. Let \[\al: G\times G\rt F^\times, \al(a,b)=(-1)^{\al_1(a,b)}. \] Then $(G,\al)$ is a symplectic abelian group.

\begin{lem}\lb{51}% (1)If \[R=\{a\in G\setminus \Rad(\al)|q(a)=1\}\] generate $G$. Then $R$ is a skew root system of Lie type in  $(G,\al)$.

 If \[R=\{a\in G|q(a)=0\}\] generate $G$. Then $R$ is a skew root system of Jordan type in  $(G,\al)$.
\end{lem}
\bp %(1) SRSL0 and SRSL1 are obvious. Now we verify SRSL2. Assume $a,b\in R$ and $\al(a,b)\neq 1$. Then $\al(a,b)=-1$, $\al_1(a,b)=1$ and $q(a+b)=q(a)+q(b)+\al_1(a,b)=1$, so $a+b\in R$.

 (SRSJ0) and (SRSJ1) are obviously satisfied by $R$. Now we verify (SRSJ2). Assume $a,b\in R$ and $\al(a,b)\neq -1$. Then $\al(a,b)=1$, $\al_1(a,b)=0$ and $q(a+b)=q(a)+q(b)+\al_1(a,b)=0$, so $a+b\in R$.
\ep

\begin{lem}\lb{52} Assume that $\al$ is nonsingular and \[R=\{a\in G\setminus \{0\}|q(a)=1\}\] generates $G$. Then $R$ is a skew root system of Lie type in the symplectic abelian group $(G,\al)$.
\end{lem}
\bp
(SRSL0) and (SRSL1) are obviously satisfied by $R$. Now we verify (SRSL2). Assume $a,b\in R$ and $\al(a,b)\neq 1$. Then $\al(a,b)=-1$, $\al_1(a,b)=1$ and $q(a+b)=q(a)+q(b)+\al_1(a,b)=1$, so $a+b\in R$.
\ep
Now assume that $n=2k+1$ is odd. One knows that there is only one type of nonsingular quadratic form on $\mathbb{F}_2^{2k+1}$ up to isomorphism, and
\be h:\mathbb{F}_2^{2k+1}\rt \mathbb{F}_2, h(a)=\sum_{i=1}^k a_{2i-1}a_{2i}+a_{2k+1}^2\ee  is such a nonsingular quadratic form. Let $\bt_1:\mathbb{F}_2^{2k+1}\times \mathbb{F}_2^{2k+1}\rt \mathbb{F}_2$ be its polarization, i.e., \[\bt_1(a,b)=h(a+b)-h(a)-h(b)=\sum_{i=1}^k (a_{2i}b_{2i-1}-a_{2i-1}b_{2i}). \] Identify the group $G=\mathbb{Z}_2^{2k+1}$ with the additive group of $\mathbb{F}_2^{2k+1}$. Let \[\bt: G\times G\rt F^\times, \bt(a,b)=(-1)^{\bt_1(a,b)}. \] One has $\Rad(\bt)=\Rad(\bt_1)=\{ 0,e_{2k+1}\}$.

\bpr\lb{57}
 Assume $k\geq 2$.

 (1) The subset \[R^-=\{a\in G\setminus \Rad(\bt)|h(a)=1\}\] is a skew root system of Lie type in $(G,\bt)$.

 (2) The corresponding reduced skew root system of $R^-$ is  $\mathbb{Z}_2^{2k}\setminus \{0\}$ in the nonsingular symplectic abelian group $\mathbb{Z}_2^{2k}$ as defined in Section 5.1. One has $R^-$ is nonreduced and indecomposable, $|R^-|=2^{2k}-1$, and $L(R^-)\cong sl(2^k,F).$
\epr
\bp
(1)(SRSL1) is clear. Now we prove (SRSL2). Note that $\Rad(\bt)=\{ 0,e_{2k+1}\}$. Assume that $a,b\in R$ and $\bt(a,b)\neq 1$, i.e., $\bt_1(a,b)\neq 0$. Then $h(a+b)=h(a)+h(b)+\bt_1(a,b)=1$. If $a+b=0$ then $a=b$ which contradicts to $\bt(a,b)\neq 1$. If $a+b=e_{2k+1}$ then $h(a)\neq h(b)$ which contradicts to $a,b\in R$. Thus $a+b\in R^-$.

Finally we prove (SRSL0). It is easy to see that the subset \[B=\{e_{2i-1}+e_{2i},\  e_{2i}+e_{2k+1}|i=1,\cdots,k\}\cup \{\sum_{j=2k-3}^{2k+1}e_j\}  \] of $R^-$ generates $G$.

Since $\Rad(\bt)\neq 0$, $R^-$ is nonreduced.

(2)
Recall $G=\mathbb{Z}_2^{2k+1}$.
As $\Rad(\bt)=\{0,e_{2k+1}\}$, $\overline{G}=G/ \Rad(\bt)= \mathbb{Z}_2^{2k}$, and the quotient group homomorphism is $p:\mathbb{Z}_2^{2k+1}\rt \mathbb{Z}_2^{2k}, (a_1,\cdots,a_{2k+1})\mapsto (a_1,\cdots,a_{2k})$. Thus  \[\overline{\bt_1}(a,b)=\sum_{i=1}^k (a_{2i}b_{2i-1}-a_{2i-1}b_{2i}) \] and \[\overline{\bt}:\overline{G}\times\overline{G}\rt F^\times, \overline{\bt}(a,b)=(-1)^{\overline{\bt_1}(a,b)} \] is nonsingular.

Next we will show that $\overline{R^-}=p(R^-)=\mathbb{Z}_2^{2k}\setminus\{0\}$. It is clear that $\overline{R^-}\subseteq \mathbb{Z}_2^{2k}\setminus\{0\}$. For any $a\in \mathbb{Z}_2^{2k} \setminus\{0\}$,
if $h(a,0)=1$ then $(a,0)\in R^-$ and $p(a,0)=a$; if $h(a,0)=0$ then $(a,1)\in R^-$ and $p(a,1)=a$. Thus $\overline{R^-}=\mathbb{Z}_2^{2k}\setminus\{0\}$, $|R^-|=|\overline{R^-}|=2^{2k}-1$, and the corresponding reduced skew root system of $(\mathbb{Z}_2^{2k+1},\bt, R^-)$ is $(\mathbb{Z}_2^{2k},\overline{\bt},\mathbb{Z}_2^{2k}\setminus\{0\})$ defined in Section 5.1. There is a surjective Lie algebra homomorphism from $L(R^-)$ to $L(\overline{R^-})$ by Lemma \ref{47}. As $|R^-|=|\overline{R^-}|$, $L(R^-)\cong L(\overline{R^-})\cong sl(2^k,F)$. Since $sl(2^k,F)$ is simple, $R^-$ must be indecomposable.
\ep
\bpr
 Assume $k\geq 1$.

(1)The subset \[R^+=\{a\in G|h(a)=0\}\] is a skew root system of Jordan type in $(G,\bt)$.

(2) The corresponding reduced skew root system of $R^+$ is $\mathbb{Z}_2^{2k}$ in the nonsingular symplectic abelian group $\mathbb{Z}_2^{2k}$ as defined in Section 5.1. One has $R^+$ is nonreduced and indecomposable, $|R^+|=2^{2k}$, and $J(R^+)$ is isomorphic to $M(2^k,F)^+$.
\epr
\bp
(1)The subset \[B=\{e_{i}|i=1,\cdots,2k\}\cup \{e_{2k-1}+e_{2k}+e_{2k+1}\} \] of $R^+$ generates $G$, thus by Lemma \ref{51} $R^+$ is a skew root system of Jordan type in $(G,\bt)$.

(2) can be proved similarly as in Proposition  \ref{57}(2).
\ep

Next we assume that $n=2k$ is even.

Let \[
\bt_1:\mathbb{F}_2^{2k}\times \mathbb{F}_2^{2k}\rightarrow \mathbb{F}_2,\bt_1(a,b)={\sum_{i=1}^k(a_{2i}b_{2i-1}-a_{2i-1}b_{2i})}
\] be the unique nonsingular alternating bilinear form on $\mathbb{F}_2^{2k}$ up to isomorphism. One knows that up to isomorphism there are two quadratic forms on $\mathbb{F}_2^{2k}$ that polarize to $\bt_1$, they are \be\lb{} f_0:\mathbb{F}_2^{2k}\rightarrow \mathbb{F}_2,f_0(a)={\sum_{i=1}^{k}a_{2i-1}a_{2i}},\ee and  \be\lb{} f_1:\mathbb{F}_2^{2k}\rightarrow \mathbb{F}_2,f_1(a)={a_1^2+a_2^2+\sum_{i=1}^{k}a_{2i-1}a_{2i}}.\ee

Let \be\lb{512}
\bt:G\times G\rightarrow F^{\times},\bt(a,b)=(-1)^{\bt_1(a,b)}.
 \ee Then $(G,\bt)$ is a nonsingular symplectic abelian group which is an orthogonal direct sum of $k$ copies of nonsingular symplectic abelian group $\mathbb{Z}_2^2$, and the twisted group algebra corresponding to  $(G,\bt)$ is isomorphic to $M(2,F)^{\ot k}\cong M(2^k,F)$.

\ble\lb{}
Assume $G=\mathbb{Z}_2^{2k}$ with $k\geq 1$. Let $R_0^-=\{a\in G\setminus \{0\}|f_0(a)=1\}$ and  $R_0^+=\{a\in G|f_0(a)=0\}$. Then

(1) If $k\geq 2$, then $R_0^-$ is a skew root system  of Lie type in $(G,\bt)$. One has $|R_0^-|=2^{2k-1}-2^{k-1}$.

(2) If $k\geq 1$, then $R_0^+$ is a skew root system  of Jordan type in $(G,\bt)$. One has $|R_0^+|=2^{2k-1}+2^{k-1}$.
 \ele
 \bp

(1) The subset \[B=\{e_{2i-1}+e_{2i},\  e_{2i-1}+e_{2i}+e_{2i+1}|i=1,\cdots,k\}  \ (e_{2k+1}=e_1)\] of $R_0^-$ generates $G$, thus by Lemma \ref{52} $R_0^-$ is a skew root system of Lie type in $(G,\bt)$.

(2) The subset \[B=\{e_{i}|i=1,\cdots,2k\} \] of $R_0^+$ generates $G$, thus by Lemma \ref{51} $R_0^+$ is a skew root system of Jordan type in $(G,\bt)$.

We omit the simple computation of the cardinality of $R_0^-$ and $R_0^+$.
\ep

\begin{lem}\lb{}
Assume  $G=\mathbb{Z}_2^{2k}$ with $k\geq 1$. Let $R_1^-=\{a\in G\setminus \{0\}|f_1(a)=1\}$ and $R_1^+=\{a\in G|f_1(a)=0\}$. Then

 (1) If $k\geq 1$, then $R_1^-$ is a skew root system of Lie type in $(G,\bt)$ and  $|R_1^-|=2^{2k-1}+2^{k-1}$.

 (2) If $k\geq 2$, then $R_1^+$ is a skew root system of Jordan type in $(G,\bt)$ and  $|R_1^+|=2^{2k-1}-2^{k-1}$.
 \end{lem}
\bp (1) It is directly verified that \bee B=\begin{cases} \{e_1\}\cup\{e_1+e_2\}, & \text{if\ k=1};\\
\{e_1,e_{2i-1}+e_{2i}|1\leq i\leq k\}\cup \{e_{2i-1}+e_{2i}+e_{2i+1}|1\leq i<k\}, & \text{if\ $k>1$.}
\end{cases} \eee is contained in $R_1^-$ and generates $G$, thus  by Lemma \ref{52} $R_1^-$ is a skew root system of Lie type in $(G,\bt)$.

(2)The subset \[B=\{e_{i}|i=3,\cdots,2k\}\cup \{e_1+e_2+e_3+e_4,e_1+e_3+e_4\} \] of $R_1^+$ generates $G$, thus by Lemma \ref{51} $R_1^+$ is a skew root system of Jordan type in $(G,\bt)$.

We omit the simple computation of the cardinality of $R_1^-$ and $R_1^+$.
\ep

Now we need to do some preparations before we can identify $L(R_i^-)$ and $J(R_i^+)$ for $i=0,1$.

 For a nonsingular symmetric or skew-symmetric bilinear form $\phi$ on $F^n$, one knows that the adjoint map $*:M(n,F)\rt M(n,F)$ defined by
  $$\phi(Xu,v)=\phi(u,X^* v) $$ is an involution (i.e., involutive anti-automorphism) on $M(n,F)$. If  $\Phi$ is the matrix of $\phi$ with respect to the standard basis of $F^n$, then in matrix form one has $X^*=\Phi^{-1} X^t \Phi.$

Conversely one can show that any involution on $M(n,F)$ can be defined in this way by some nonsingular symmetric or skew-symmetric bilinear form on $F^n$. See Section 5 of \cite{bsz} for details. We call $(M(n,F),*)$ an \textit{involutive} matrix algebra. In the remaining of the section, $M$ will be assumed to be a matrix algebra.

Assume that $(M,*)$ is an {involutive} matrix algebra. Let \[\Aut(M,*)=\{\si\in \Aut(M)|\si\circ *=*\circ\si\}.\] If $H$ is a subgroup of $\Aut(M,*)$, then its action on $M$ is compatible with $*$,i.e.,$h\circ *=*\circ h$ for any $h\in H$, and we say that $(M,*)$ has an $H$-action. If $M$ has a $G$-grading compatible with $*$,i.e., $M_g^*=M_g$ for any $g\in G$, then we say that $(M,*)$ has a $G$-grading. It is clear that $(M,*)$ has a $G$-grading if and only if $(M,*)$ has a $\wh{G}$-action.

Let \[K(M,*)=\{X\in M|X^*=-X\}\] and \[H(M,*)=\{X\in M|X^*=X\}.\] The subspace $K(M,*)$ is closed under the Lie bracket $[a,b]=ab-ba$ and is always regarded as a Lie algebra with this Lie bracket.  The subspace $H(M,*)$ is closed under the Jordan bracket $a\circ b=\frac{1}{2}(ab+ba)$ and is always regarded as a Jordan algebra in this way. A $G$-grading on $(M,*)$ clearly induces a $G$-grading on $K(M,*)$ (resp. on $H(M,*)$).

If $*$ corresponds to some nonsingular symmetric bilinear form on $F^n$, then $K(M(n,F),*)$ is the simple Lie algebra isomorphic to $so_n(F)$ (which is simple if $n=3$ or $n\geq 5$) and $H(M(n,F),*)$ is a Jordan algebra of type $H(F_n)$ (which is simple if $n\geq 3$).
If $*$ corresponds to some nonsingular skew-symmetric bilinear form on $F^{2n}$, then $K(M({2n},F),*)$ is the simple Lie algebra isomorphic to $sp_{2n}(F)$ and $H(M(n,F),*)$ is a Jordan algebra of type $H(Q_{2n})$ (which is simple if $n\geq 3$).

 Let $G=\mathbb{Z}_2\times \mathbb{Z}_2$,$M=M(2,F)$. Let
\be X=\left(
             \begin{array}{cc}
               -1 & 0\\
               0 & 1 \\
             \end{array}
           \right) \an \
           Y=\left(
             \begin{array}{cc}
               0 & 1\\
               1 & 0 \\
             \end{array}
             \right).
\ee
 %Let $Z_{(i,j)}=X_2^i Y_2^j$ for $(i,j)\in G$.
 By (\ref{54}), for suitable $\xi\in L^2(G,F^\times)$ one has \be \lb{56}\eta: F^\xi (\mathbb{Z}_2\times \mathbb{Z}_2)\rt M(2,F), \overline{(i,j)}\mapsto X^i Y^j\ee is a graded algebra isomorphism.
%(The grading \[M_2(F)=\op_{(i,j)\in G} F\cdot X^i Y^j\] is referred as the standard $G$-grading on $M_2(F)$.)

Assume that $*$ is an involution on $M$ corresponding to the skew-symmetric form with the  matrix $\Phi_1=\left(
             \begin{array}{cc}
               0 & 1\\
               -1 & 0 \\
             \end{array}
           \right)$. By Lemma 3 (1) in \cite{bsz}, $(M,*)$ has the $G$-grading with \be \lb{a1}K(M,*)=\{\left(
             \begin{array}{cc}
               a & b\\
               c & -a \\
             \end{array}
           \right)\}=\bigoplus_{(i,j)\in \mathbb{Z}_2^2\setminus (0,0)} F\cdot X^i Y^j \ee and \[H(M,*)=\{\left(
             \begin{array}{cc}
               a & 0\\
               0 & a \\
             \end{array}
           \right)\}= F\cdot X^0 Y^0.\]

 Assume that $\tilde{*}$ is an involution on $M$ corresponding to the symmetric form with the matrix $I_2=\left(
             \begin{array}{cc}
               1 & 0\\
               0 & 1 \\
             \end{array}
           \right)$.

By Lemma 3 (3) of \cite{bsz}, $(M,\tilde{*})$ has the $G$-grading with \be \lb{a2}K(M,\tilde{*})=\{\left(
             \begin{array}{cc}
               0 & b\\
               -b & 0 \\
             \end{array}
           \right)\}= F\cdot X^1 Y^1\ee and \[H(M,\tilde{*})=\{\left(
             \begin{array}{cc}
               a & b\\
               b & c \\
             \end{array}
           \right)\}=\bigoplus_{(i,j)\in \mathbb{Z}_2^2\setminus (1,1)} F\cdot X^i Y^j.\]
\bigskip

Next we introduce the direct product gradings on tensor product of associative algebras. Assume for $i=1,2$, $A_i$ is an associative algebra, and $(\Ga_i,G_i)$ is the respective gradings on $A_i$, then $(\Ga, G_1\times G_2)$ is the \textit{direct product grading} on $A=A_1\otimes A_2$ with $$\Gamma:A_1\ot A_2=\bigoplus_{(a,b)\in G_1\times G_2 } (A_1\otimes A_2)_{(a,b)}$$ where $$(A_1\otimes A_2)_{(a,b)}=(A_1)_a\otimes (A_2)_b, (a,b)\in G_1\times G_2.$$

If $(M_i,*_i)$ is an involutive matrix algebra for $i=1,2$, then their direct product $(M_1,*_1)\ot (M_2,*_2)$ is defined to be the involutive matrix algebra $(M_1\ot M_2,*_1\ot *_2)$, where $(C\ot B)^{*_1\ot *_2}=C^{*_1}\ot B^{*_2}$ for any $C\ot B\in M_1\ot M_2$. If $(M_i,*_i)$ has a $G_i$-grading, then the $G_1\times G_2$-grading on $M_1\ot M_2$ is compatible with $*_1\ot *_2$, thus $(M_1\ot M_2,*_1\ot *_2)$ has the naturally defined $G_1\times G_2$-grading.

Assume that $(M_i,*_i)$ has a $\mathbb{Z}_2\times \mathbb{Z}_2$-grading for $i=1,\cdots,k$, where $M_i=M(2,F)$ and the $\mathbb{Z}_2\times \mathbb{Z}_2$-grading on $M_i$ is as in (\ref{56}) for all $i$. Then the tensor product $(M,*)$ of $(M_i,*_i)$ for $i=1,\cdots,k$ has a $G$-grading with \[M\cong M(2^k,F), \ G=(\mathbb{Z}_2\times \mathbb{Z}_2)^k=\mathbb{Z}_2^{2k}.\]

By (\ref{56}) and Lemma \ref{m}, one has the $G$-grading isomorphism: \be\lb{55} \psi: F^\xi G\rt M(2^k,F), \overline{(i_1,j_1,\cdots,i_k,j_k)}\mapsto X^{i_1}Y^{j_1}\ot\cdots\ot X^{i_k}Y^{j_k}.\ee

Since the grading on $M$ is compatible with the involution $*$, $K(M,*)$ (resp. $H(M,*)$) also has a $G$-grading. Assume that $*_i$ corresponds to a nonsingular bilinear form $\phi_i$ on $F^2$, then $*$ corresponds to the nonsingular bilinear form $\phi=\phi_1\ot\cdots\ot\phi_k$ on $(F^2)^{\ot k}\cong F^{2^k}$.  Assume the total number of skew-symmetric factors in $\phi$ is $m$. If $m$ is even, then $\phi$ is symmetric, one has $K(M,*)\cong so(2^k,F)$ and $H(M,*)$ is a Jordan algebra of type $H(F_{2^k})$. If $m$ is odd, then $\phi$ is  skew-symmetric, one has $K(M,*)\cong sp(2^k,F)$ and $H(M,*)$ is a Jordan algebra of type $H(Q_{2^k})$. Next we will compute the support $R^-$ of the $G$-grading on $K(M,*)$ and the support $R^+$ of the $G$-grading on $H(M,*)$ in the two cases, $m=0$ or $m=1$, and get the desired results. Since the $G$-grading on $M$ is multiplicity free, $R^-$ and $R^+$ are disjoint and their union is $G$.

%Recall the nonsingular quadratic forms $g$ and $f$ on $\mathbb{F}_2^{2k}$ defined in (\ref{aa8}) and (\ref{aa9}). One knows that they both polarize to $\bt_1$.

 %Let  \be\lb{l5} \bt: \mathbb{F}_2^{2k}\times  \mathbb{F}_2^{2k}\rightarrow F^{\times},\bt(a,b)=(-1)^{\bt_1(a,b)}=(-1)^{{\sum_{i=1}^k(a_{2i}b_{2i-1}-a_{2i-1}b_{2i})}}.\ee
\bpr
Let $G=\mathbb{Z}_2^{2k}$ and $\bt$ be the nonsingular symplectic bilinear form on $G$ as in (\ref{512}).

(1) Assume  the total number of skew-symmetric factors $\phi_i$ in $\phi$ is 0 and the matrix for each $\phi_i$ is $I_2$. Then $K(M,*)\cong so(2^k,F)$ and $H(M,*)$ is a Jordan algebra of type $H(F_{2^k})$.

(1.1)The support of $K(M,*)$ is $R_0^-=\{a\in G\setminus \{0\}|f_0(a)=1\}$ and the support of $H(M,*)$ is $R_0^+=\{a\in G|f_0(a)=0\}$.

(1.2)Conversely, one knows that $R_0^-$ is a skew root system of Lie type in $(G,\bt)$  for $k\geq 2$ and $R_0^+$ is a skew root system of Jordan type in $(G,\bt)$ for $k\geq 1$. The $G$-grading isomorphism $\psi$ in (\ref{55}) induces the Lie algebra isomorphism $L(R_0^-)\rt K(M,*)\cong so(2^k,F)$ for $k\geq 2$ and the Jordan algebra isomorphism $J(R_0^+)\rt H(M,*)$ for $k\geq 1$.

(2) Assume the total number of skew-symmetric factors $\phi_i$ in $\phi$ is 1, the matrix of $\phi_1$ is $\Phi_1$ and the matrix of $\phi_i$ is  $I_2$ for $i>1$. Then $K(M,*)\cong sp(2^k,F)$ and $H(M,*)$ is a Jordan algebra of type $H(Q_{2^k}$).

(2.1) The support of $K(M,*)$ is $R_1^-=\{a\in G\setminus \{0\}|f_1(a)=1\}$ and the support of $H(M,*)$ is $R_1^+=\{a\in G|f_1(a)=0\}$.

(2.2) Conversely, one knows that $R_1^-$ is a skew root system of Lie type in $(G,\bt)$  for $k\geq 1$ and $R_1^+$ is a skew root system of Jordan type in $(G,\bt)$ for $k\geq 2$. The $G$-grading isomorphism $\psi$ in (\ref{55}) induces the Lie algebra isomorphism $L(R_1^-)\rt K(M,*)\cong sp(2^k,F)$ for $k\geq 1$ and the Jordan algebra isomorphism $J(R_1^+)\rt H(M,*)$ for $k\geq 2$.
\epr
\bp
Since (1.2) (resp. (2.2)) follows from (1.1) (resp. (2.1)) by Lemma \ref{461}, we only need to prove (1.1) and (2.1).

First let us consider the case $k=1$. Then $G=\mathbb{Z}_2\times \mathbb{Z}_2$ and $M=M(2,F)$. Then by (\ref{a2}) the support of $K(M,\tilde{*})=so(2,F)$ is \be\lb{m2}\{(1,1)\}=\{x\in G|f_0(x)=x_1x_2=1\}.\ee Then the support of  $H(M,\tilde{*})$ is
$\{x\in G|f_0(x)=0\}$.

By (\ref{a1}) the support of $K(M,*)=sp(2,F)$ is \be\lb{m3}\{(0,1),(1,0),(1,1)\}=\{x\in G|f_1(x)=x_1^2+x_2^2+x_1x_2=1\}.\ee Then the support of  $H(M,*)$ is
$\{x\in G|f_1(x)=0\}$.

Next let us consider the general case. Then $G=\mathbb{Z}_2^{2k}$. Let $C=C_1\ot\cdots\ot C_k\in M(2,F)\ot\cdots\ot M(2,F)$, $C_i$ being one of $X^i Y^j$ with $(i,j)\in \mathbb{Z}_2\times \mathbb{Z}_2$. Let $v_1\ot\cdots\ot v_k, u_1\ot\cdots\ot u_k\in (\mathbb{F}^2)^{\ot k}$. Then \be\lb{m1}\phi_i(C_i(v_i),u_i)=\begin{cases} -\phi_i(v_i,C_i(u_i)), &\text{if $C_i\in K(M(2,F), *_i)$;}\\ \phi_i(v_i,C_i(u_i)), &\text{if $C_i\in H(M(2,F), *_i)$.}\end{cases}\ee
One has
 \[\phi( C_1\ot\cdots\ot C_k(v_1\ot\cdots\ot v_k),u_1\ot\cdots\ot u_k)=\phi_1(C_1(v_1),u_1)\cdots \phi_k(C_k(v_k),u_k),\] and
\[\phi(v_1\ot\cdots\ot v_k, C_1\ot\cdots\ot C_k(u_1\ot\cdots\ot u_k))=\phi_1(v_1,C_1(u_1))\cdots \phi_k(v_k,C_k(u_k)).\]
Thus by (\ref{m1}) \[\phi( C_1\ot\cdots\ot C_k(v_1\ot\cdots\ot v_k),u_1\ot\cdots\ot u_k)\]
\[=(-1)^t \phi(v_1\ot\cdots\ot v_k, C_1\ot\cdots\ot C_k(u_1\ot\cdots\ot u_k))\]
where $t$ is the number of $C_i\in K(M, *_i)$.

(1) Assume  the total number of skew-symmetric factors $\phi_i$ in $\phi$ is 0 and the matrix for each $\phi_i$ is $I_2$.  Recall \[f_0:\mathbb{F}_2^{2k}\rt \mathbb{F}_2,f_0(x)=\sum_{i=1}^k x_{2i-1}x_{2i}.\] Let $a=(a_1,a_2,\cdots,a_{2k-1},a_{2k})\in \mathbb{Z}_2^{2k}$. Then by (\ref{m2}) $a$ is in the support of $K(M,*)$ if and only if $f_0(a)=1$ and $a$ is in the support of $H(M,*)$ if and only if $f_0(a)=0$. So the support of $K(M,*)$ is those $a\in G$ with $f_0(a)=1$ and the support of $H(M,*)$ is those $a\in G$ with $f_0(a)=0$.

(2) Assume the total number of skew-symmetric factors $\phi_i$ in $\phi$ is 1, the matrix of $\phi_1$ is $\Phi_1$ and the matrix of $\phi_i$ is  $I_2$ for $i>1$. Recall \[f_1:\mathbb{F}_2^{2k}\rt \mathbb{F}_2,f_1(x)=(x_1^2+x_2^2+x_1 x_2)+\sum_{i=2}^k x_{2i-1}x_{2i}=x_1^2+x_2^2+\sum_{i=1}^k x_{2i-1}x_{2i}.\]  Let $a=(a_1,a_2,\cdots,a_{2k-1},a_{2k})\in G$. Then  by (\ref{m2}) and (\ref{m3}) $a$ is in the support of $K(M,*)$ if and only if $f_1(a)=1$ and $a$ is in the support of $H(M,*)$ if and only if $f_1(a)=0$. So the support of $K(M,*)$ is those $a\in G$ with $f_1(a)=1$ and the support of $H(M,*)$ is those $a\in G$ with $f_1(a)=0$.

\ep

We have constructed 3 families of skew root systems of Lie type (resp. of Jordan type), and identified the corresponding Lie algebras $L(R)$ (resp. the corresponding Jordan algebras $J(R)$) with the multiplicity free semisimple grading (resp. with the multiplicity free division grading). As far as we know, all the multiplicity free semisimple gradings on simple Lie algebras and multiplicity free division gradings on special simple Jordan algebras are included in them. In particular, exceptional simple Lie algebras do not have multiplicity free semisimple gradings. We believe that all the multiplicity free semisimple gradings on simple Lie algebras and all the multiplicity free division gradings on special simple Jordan algebras can be constructed by means of skew root system of Lie type and of Jordan type respectively.

\end{document}